# On Cantor's Theorem

W. Mückenheim

University of Applied Sciences, Baumgartnerstraße 16, D-86161 Augsburg, Germany

mueckenh@rz.fh-augsburg.de

________________________________________

**Abstract.** The famous contradiction of a bijection between a set and its power set is a consequence of the impredicative definition involved. This is shown by the fact that a simple mapping between equivalent sets does also fail to satisfy the critical requirement.

## 1. Introduction

There is no surjection of the set $\mathbb{N}$ of all natural numbers on its power set $\mathcal{P}(\mathbb{N})$. This proof was first given by Hessenberg [1]. Though Zermelo [2], calling it Cantor's theorem, attributed it to Cantor [3] himself, Cantor's most famous paper [3] does not contain the notion power set (Potenzmenge) at all. But the expression "Cantor's theorem" has become generally accepted.

If there was a surjective mapping of $\mathbb{N}$ on its power set, $s: \mathbb{N} \to \mathcal{P}(\mathbb{N})$, then some natural numbers $n \in \mathbb{N}$ might be mapped on subsets $s(n)$ not containing $n$. Call these numbers $n$ "non-generators" to have a convenient abbreviation. The subset $M$ of all non-generators

$$M = \{n \in \mathbb{N} \mid n \notin s(n)\}$$

belongs to $\mathcal{P}(\mathbb{N})$ as an element, whether containing numbers or being empty. But there is no $m \in \mathbb{N}$ available to be mapped on $M$. If $m$ is not in $M$, then $m$ is a non-generator, but then $M$ must contain $m$ and vice versa: $m \in M \Rightarrow m \notin M \Rightarrow m \in M \Rightarrow \ldots$ This condition, however, has nothing to do with the cardinal numbers of $\mathbb{N}$ and $\mathcal{P}(\mathbb{N})$ but it is simply a paradox request: *m has to be mapped by s on a set which does not contain it, if it contains it*.

The set {$M$, $m$, $s$} simply belongs to the class of impredicatively defined sets like Russell's set of all sets which do not contain themselves. Those types of definition are well known for their power of generating paradoxes and have been banned from set theory long ago. The set of all ordinals, the set of all cardinals, and the set of all sets are such "impossible sets". Only the set {$M$, $m$, $s$} has survived.

## 2. A surjective mapping with a paradox condition

We will see now, that the impossible set does not exist and that the paradox-generating requirement cannot be satisfied, even if the mapping is defined between equivalent sets.

Define a bijective mapping from {1, $a$} on $\mathcal{P}(\{1\}) = \{\{\}, \{1\}\}$, where $a$ is a symbol but not a number. There are merely two bijections possible. The set $M$ of all *numbers* which are non-generators cannot be mapped by a *number m* although $M$ is in the image of both the possible mappings.

$f$: 1 → {1} and $a$ → { }  with $M_f$ = { },

$g$: 1 → { } and $a$ → {1}  with $M_g$ = {1}.

Here we have certainly no problem with lacking elements in the domain. Nevertheless Hessenberg's condition cannot be satisfied. Both the sets {$M_f$, $m_f$, $f$} and {$M_g$, $m_g$, $g$} with $m = 1$, the only available *number*, are impossible sets. Hessenberg's proof does not concern the question whether or not $\aleph_0 < 2^{\aleph_0}$.